\documentclass[a4paper,12pt]{article}
\usepackage{amsmath}     
\usepackage{amsthm}     
 
\newcommand{\beq}{\begin{equation}}
\newcommand{\eeq}{\end{equation}}

\newcommand{\beqnn}{\begin{equation*}}
\newcommand{\eeqnn}{\end{equation*}}
\DeclareMathOperator*{\mgrlongrightarrow}{\longrightarrow}

\title{Metric on state space of Markov chain}
\author{M.G.~Rozinas \footnote{e-mail: rozinas@gmail.com}}
\date{}

\begin{document}
\maketitle

\begin{abstract}
We consider finite irreducible Markov chains. 
It was shown that mean hitting time from one state to another satisfies the triangle inequality. 
Hence, sum of mean hitting time between couple of states in both directions is a metric on the space of states.
\end{abstract}

\newtheorem{corollary}{Corollary}
\newtheorem{lemma}{Lemma}
\newtheorem{theorem}{Theorem}


\section{Paths, sets of paths and probabilities}
Basic definitions and all used properties of Markov chains could be found in \cite{cite_1}. 

Consider Markov chain $\mathbf{M}$ with finite state space $Z=\{z_1,\ldots z_n\}$   
and matrix of transition probabilities \mbox{$P=(p_{i,j})_{i,j=1,\ldots,n}$}; $P$ is non-negative and for all $i=1,\ldots,n$  
    
\beq \label{Formula_01}
	\sum_{j=1}^n p_{i,j}=1  
\eeq

Let $\lambda^m=(\lambda_1^m,\ldots,\lambda_n^m)$ is distribution of probabilities of states on step $m$ ($\lambda^0$ is initial distribution). It is known that for all $m\geq0$ 

\beqnn
	\lambda^m=\lambda^0 P^m
\eeqnn

Let  $P^m=(p_{i,j}^{(m)})_{i,j=1,\ldots,n}$. Coefficient $p_{i,j}^{(m)}$  is the probability that chain with initial state $z_i$ will be in the state $z_j$ after $m$ steps.

Any sequence  $x=(x_0,\ldots\,x_m)$ of states we call a \emph{path} of length $m$. Denote length of path and its first and last elements as $Len(x)=m$, $Begin(x)=x_0$ and $End(x)=x_m$.

If $End(x)=Begin(y)$, then $x \otimes y$ is \emph{concatenation of paths $x$ and $y$}. 

Path $y$ is \emph{extension} of path $x$, if 
\beqnn
	(Len(x) \leq Len(y)) \quad  \& \quad (\forall i \leq Len(x))(x_i=y_i) 
\eeqnn
 
Let $S$ is the set of all paths and, for any $m \geq 0$, $a,b \in Z$, 

\beqnn
	S_{(a,b)}^m =  \{s\in S | (Len(s)=m) \& (Begin(s)=a)  \&  (End(s)=b)\}
\eeqnn

Asterisk instead of some index in $S_{(a,b)}^m$ means dropping of corresponding condition. For example, 
$S_{(a,*)}^* = \{s\in S \mid Begin(s)=a\}$ .

For any path $x$ denote by $Ext^k(x)$ set of all extensions of path $x$, having length $k$:
\beqnn
	Ext^k(x)  =  \{ s \in S  | (Len(s)=k)  \&  (s \text{ is extension of } x) \}
\eeqnn
 
For any set of paths $X \subseteq S $, let
\beqnn
	Ext^k(X)= \bigcup_{x\in X} Ext^k(x)	
\eeqnn                                        
Sometimes we will write $p(z_i,z_j)$ instead of $p_{i,j}$.  	
For any $a \in Z$, $m>0$, $x \in S_{(a,*)}^m$ define 
\beq \label{Formula_02_5}
	P_a^m(x) = \prod_{i=1}^m p(x_{i-1},x_i)
\eeq
Evidently, $P_a^m(x)$ is the probability that Markov chain $\mathbf{M}$, being at initial moment in state $a$, will follow path $x$ on first $m$ steps.

Similarly, for any $X \subseteq S_{(a,*)}^m$
\beqnn
	P_a^m(X) = \sum_{x \in X} P_a^m(x)
\eeqnn
is the probability that Markov chain $\mathbf{M}$, being at initial moment in state $a$, will follow on first $m$ steps to some path from $X$. Note that 
\beqnn
	p_{a,b}^{(m)} = P_a^m(S_{(a,b)}^m)
\eeqnn
and 
\beq \label{Formula_04}
	P_a^m(S_{(a,*)}^m) = \sum_{b \in Z} P_a^m(S_{(a,b)}^m) = 1 
\eeq

For any $X \subseteq S$ define function $P$ on sets of paths $X$ as 

\beqnn
	P(X)= \sum_{x \in X} P(x) =  \sum_{x \in X} P_{Begin(x)}^{Len(x)}(x)  
\eeqnn

\begin{lemma} \label{lemma:L01}
	If $X \in S_{(a,*)}^m$ and $k\geq m$, then $P_a^k(Ext^k(X))=P_a^m(X)$
\end{lemma}
 
\begin{proof}

Case $k=m$ is trivial.  For $k=m+1$   

\beqnn
 	P_a^{m+1} (Ext^{m+1}(X))  = \sum_{y \in Ext^{m+1}(X)} P_a^{m+1}(y) = 
 	\sum_{x \in X , z\in Z} P_a^{m+1}(x\otimes (x_m,z))=  
\eeqnn

\beqnn
   \sum_{x \in X} \sum_{z \in Z}(( \prod_{i=1}^m p(x_{i-1},x_i))p(x_m,z))=  \sum_{x \in X}( ( \prod_{i=1}^m p(x_{i-1},x_i)) \sum_{z \in Z} p(x_m,z))   
\eeqnn
  
According to \eqref{Formula_01}, $$\sum_{z \in Z} p(x_m,z) = 1$$ therefore 
\beqnn
 	P_a^{m+1} (Ext^{m+1}(X))  = \sum_{x \in X} \prod_{i=1}^m p(x_{i-1},x_i)= \sum_{x \in X} P_a^m(x) = P_a^m(X)
\eeqnn
 
For any $k \geq m$ extension  $Ext^k(X)$ could be obtained by   successive one-step extensions, therefore lemma hold for all $k \geq m$.

\end{proof}


\section{Irreducible finite Markov chains}

Markov chain $\mathbf{M}$ is \emph{irreducible}, if there is positive probability of transition from any state to any other (possibly, in more than one step), i.e.

\beqnn
	(\forall i,j)(\exists m)(p_{i,j}^{(m)}>0) 
\eeqnn

For any set of paths $X$ and any $W \subseteq Z$ denote by $X_{[-W]}$ the set of paths from $X$ that does not contain states from $W$, and by $X_{[+W]}$ the set of paths from $X$ that contain at least one element from $W$:
   
\beqnn
	X_{[-W]}=\{ x  \in X \mid (\forall w \in W)(w \notin X) \}
\eeqnn

\beqnn
	X_{[+W]}=\{ x  \in X \mid (\exists w \in W)(w \in X) \}
\eeqnn
 
For one-element sets we will use shortcuts $X_{[\pm w]}=X_{[\pm \{ w\} ]}$.

\begin{lemma} \label{lemma:L02}   
	If finite Markov chain is irreducible, then for any $a,w \in Z$ 	
\beqnn
	\lim_{m \to \infty} P_a^m (S_{(a,*)[-w]}^m  )=0
\eeqnn
\end{lemma}

\begin{proof}

For any $u,v \in Z$ expression $P_a^m (S_{(u,*)[+v]}^m)$ is monotonously non-strictly increasing function of $m$ or, equivalently, 
$P_a^m (S_{(u,*)[-v]}^m)$ is monotonously non-strictly decreasing. Actually, 
$$
m_1<m_2 \quad \Rightarrow \quad
Ext^{m_2}(S_{(u,*)[+v]}^{m_1}) \subseteq S_{(u,*)[+v]}^{m_2} \quad \Rightarrow \quad
P_u^{m_2}(Ext^{m_2}(S_{(u,*)[+v]}^{m_1}))   \leq P_u^{m_2}(S_{(u,*)[+v]}^{m_2})
$$ 
According to Lemma \ref{lemma:L01},
$P_u^{m_1}(S_{(u,*)[+v]}^{m_1})) = P_u^{m_2}(Ext^{m_2}(  S_{(u,*)[+v]}^{m_1}))$,
therefore  
\beqnn
	P_u^{m_1}(S_{(u,*)[+v]}^{m_1})) \leq P_u^{m_2}(S_{(u,*)[+v]}^{m_2}))
\eeqnn

Due to irreducibility, for all $u,v \in Z$ exists index $t(u,v)$ that satisfies condition $p_{u,v}^{(t(u,v))} > 0 $.
Let $t=max\{ t(u,v) \}$ and $q = min\{ p_{u,v}^{(t(u,v))} \} $ . Then, for any $u,v \in Z$, 
\beqnn
	q \leq p_u^{(t(u,v))} = P_u^{t(u,v)}(S_{(u,v)}^{t(u,v)}) \leq P_u^{t(u,v)}(S_{(u,*)[+v]}^{t(u,v)}) \leq  P_u^t(S_{(u,*)[+v]}^t)
\eeqnn
Hence, for some $q>0$ and any $u,w \in Z$  
\beq \label{Formula_05}
	P_u^t(S_{(u,*)[-w]}^t) \leq 1-q
\eeq
Inequality \eqref{Formula_05} means that for any sequence of $t$ steps, probability of event ``all states in the
sequence differ from $w$'' does not exceed $1-q$. Hence, for any sequence, containing $kt$  steps, probability of avoiding state $w$ does not exceed $(1-q)^k$, and for any $a,w \in Z$
\beqnn 
	P_a^m(S_{(a,*)[-w]}^m) \leq (1-q)^{\left\lfloor m/t \right\rfloor} \mgrlongrightarrow_{m \to \infty}0 
\eeqnn
\end{proof}

\begin{corollary} \label{Cor01} 
If finite Markov chain is irreducible, then exist $t_0$, $\alpha>0$, $1>\beta>0$ such that for any $u,v \in Z$ and any $t \geq t_0$  
\beqnn
	P_u^t(S_{(u,*)[-v]}^t)\leq\beta^{\alpha t}
\eeqnn
\end{corollary}

\begin{proof}
Proof is evident from the proof of lemma \ref{lemma:L02}.
\end{proof}

Let $X$ is a set of paths and $b \in W \subseteq Z$. Denote by 
$X_{[+W:b]}$ the set of paths from X that contain elements from $W$, and $b$ is the first element on the path that belongs to $W$:
  
\beqnn
	X_{[+W:b]}=\{ x \in X_{[+W]} \mid (\forall c \in W)(x \in X_{[+c]} \longrightarrow x \in X_{[+b]}) \}
\eeqnn

\begin{lemma} \label{lemma:L04} 
If finite Markov chain is irreducible, $\oslash \subset W \subseteq Z$,  and $a \in Z$, then for any $b \in W$ exist limit
\beqnn
  \overline{p}_{{}_W,a,b}=\lim_{m \to \infty}P_a^m(S_{(a,*)[+W:b]}^m) 
\eeqnn
and
\beq \label{Formula_05.5}
	\sum_{b \in W}{\overline{p}_{{}_W,a,b}}=1	
\eeq
\end{lemma}

\begin{proof}																 
If $a \in W$, statement is trivial and 
 
\beqnn
\overline{p}_{{}_W,a,b} =
\begin{cases}
1 & \text{if } b=a\\
0 & \text{if } b\neq a
\end{cases}
\eeqnn

For any $a \notin W $, $b \in W$, expression
$P_a^m(S_{(a,*)[+W:b]}^m)$  is monotonously non-strictly increasing function of $m$. Really, if $m_1<m_2$, then
$Ext^{m_2}(S_{(a,*)[+W:b]}^{m_1}) \subseteq S_{(a,*)[+W:b]}^{m_2} $  
and, according to Lemma \ref{lemma:L01}, 
$ P_a^{m_1}(S_{(a,*)[+W:b]}^{m_1}) = P_a^{m_2}(Ext^{m_2}(S_{(a,*)[+W:b]}^{m_1})) \leq P_a^{m_2}(S_{(a,*)[+W:b]}^{m_2})$.

Evidently, sets $\{S_{(a,*)[+W:b]}^m\}_{b \in W}$ are pairwise disjoint subsets of $S_{(a,*)}^m$  and hence, according to \eqref{Formula_04},
\beq \label{Formula_06}
	\sum_{b \in W}{P_a^m(S_{(a,*)[+W:b]}^m)} \leq 1
\eeq 	  
Because addends in \eqref{Formula_06} are non-negative monotonously non-decreasing functions of $m$, limits 
  $$\overline{p}_{{}_W,a,b}=\lim_{m \to \infty}P_a^m(S_{(a,*)[+W:b]}^m)$$ 
exist and 
\beqnn 
	\sum_{b \in W}{\overline{p}_{{}_W,a,b}} \leq 1	
\eeqnn
Note that 
$$ S_{(a,*)}^m = \big( \bigcup_{b \in W}{S_{(a,*)[+W:b]}^m} \big) \cup S_{(a,*)[-W]}^m $$ 
therefore 
\beq \label{Formula_06.5}
 \sum_{b \in W}{P_a^m(S_{(a,*)[+W:b]}^m)}
=1-P_a^m(S_{(a,*)[-W]}^m  ) 
\eeq

Select some element $c$ in non-empty set $W$. From $c \in W$ follows $S_{(a,*)[-W]}^m \subseteq  S_{(a,*)[-c]}^m $
and  $P_a^m(S_{(a,*)[-W]}^m) \leq  P_a^m(S_{(a,*)[-c]}^m)$. 

According to lemma \ref{lemma:L02}, 
$\lim_{m \to \infty}{P_a^m(S_{(a,*)[-c]}^m)}=0$, therefore 

\beq \label{Formula_06.6}
	\lim_{m \to \infty}{P_a^m(S_{(a,*)[-W]}^m)}=0
\eeq
From \eqref{Formula_06.5} and \eqref{Formula_06.6} follows 
\beqnn 
	\lim_{m \to \infty}{\sum_{b \in W}{P_a^m(  S_{(a,*)[+W:b]}^m)}}=1 
\eeqnn
and \eqref{Formula_05.5} proved.
\end{proof}

\begin{corollary} \label{Cor02} 
If finite Markov chain is irreducible, $\oslash \subset W \subseteq Z$,  and $a \in Z$, then for any $b \in W$ and any $m \geq 1$ 
\beqnn
  \overline{p}_{{}_W,a,b} \geq P_a^m(S_{(a,*)[+W:b]}^m) 
\eeqnn
\end{corollary}

\begin{proof}
	Proof is evident from the proof of lemma \ref{lemma:L04}.
\end{proof}

Let $W \subseteq Z$, $ |W|>1$ , $a,b \in W$ and $a \neq b$.  
We say that path $x$ is \emph{$W$-arrow from $a$ to $b$}, if 
\beqnn
	(Begin(x)=a) \quad \& \quad (End(x)=b) \quad \& \quad (\forall i<Len(x))(x_i \in W \Rightarrow x_i=a)
\eeqnn
i.e. path $x$ leads from $a$ till the first hitting $b$, and contains no states from $W$ besides $a$ and $b$.  
Denote by $R_{a,b}^W$ the set of all $W$-arrows from $a$ to $b$.

\begin{lemma} \label{lemma:L06} 
If finite Markov chain is irreducible, $W \subseteq Z$, $|W|>1$, 
$a,b \in Z$ and $a \neq b$, then
\beqnn
  \overline{p}_{{}_{W \backslash \{a\},a,b}} = P(R_{a,b}^W) 
\eeqnn
\end{lemma}

\begin{proof}

By definition, 
\beqnn
	\overline{p}_{{}_{W \backslash \{a\},a,b}} = \lim_{m \to \infty}{P_a^m(S_{(a,*)[+ W \backslash \{a\} : b]}^m)}
\eeqnn
Note that $S_{(a,*)[+ W \backslash \{a\} : b]}^m$ is the set of paths of length $m$ that start in $a$, hit some others states from  
$W \backslash \{a\}$ ,  and $b$ is the first one among them. This means that 
\beqnn
	S_{(a,*)[+ W \backslash \{a\} : b]}^m = \bigcup_{(x \in  R_{a,b}^W)\&(Len(x) \leq m)}{Ext^m(x)}
\eeqnn 
All sets of paths $Ext^m(x)$ in the union on the right part are pairwise disjoint, therefore according to lemma \ref{lemma:L01},  
$$
P_a^m(S_{(a,*)[+ W \backslash \{a\} : b]}^m) = 
P_a^m(\bigcup_{{}_{Len(x) \leq m}^{x \in  R_{a,b}^W}}{Ext^m(x)})=
\sum_{{}_{Len(x) \leq m}^{x \in  R_{a,b}^W}}{P_a^m(Ext^m(x))}=
\sum_{{}_{Len(x) \leq m}^{x \in  R_{a,b}^W}}{P_a^{Len(x)}(x)}
$$
Hence
\beqnn
	\lim_{m \to \infty}{P_a^m(S_{(a,*)[+ W \backslash \{a\} : b]}^m)}=
	\lim_{m \to \infty}{\sum_{{}_{Len(x) \leq m}^{x \in  R_{a,b}^W}}{P_a^{Len(x)}(x)}}=
\sum_{x \in  R_{a,b}^W}{P_a^{Len(x)}(x)}=P(R_{a,b}^W)
\eeqnn
and proof completed.
\end{proof}

According to lemmas \ref{lemma:L02},\ref{lemma:L04},\ref{lemma:L06},
$P(R_{a,b}^W)$ is conditional probability of event 
``$b$ will be the first touched state from $W \backslash \{a\}$,
under condition that initial state is $a$''.


\section{Factor chain}

Using lemma \ref{lemma:L04}, introduce, for finite irreducible Markov chains, notion of \emph{factor chain}. Let $W \subseteq Z$ and $|W|>1$.  Consider new Markov chain $\mathbf{\overline{M}}$ with the state space $W$ and transition matrix 
$\overline{P}=(\overline{p})_{a,b \in W}$, 
where   
\beq \label{Formula_07}
	\overline{p}_{a,b}=\begin{cases}
P(R_{a,b}^W) & \text{if } a\neq b\\
0 & \text{if } a = b
\end{cases}
\eeq

We will say that $\mathbf{\overline{M}}$ is \emph{factor chain of $\mathbf{M}$ by set of states $W$} and write 
\mbox{$\mathbf{\overline{M}} = \mathbf{M}/W$}.
Intuitively, moving from $\mathbf{M}$ to 
$\mathbf{\overline{M}} = \mathbf{M}/W$ means

\begin{enumerate}
\item[-] Ignore all states that not belong to $W$,
\item[-] Consider as single step any sequence of steps from some element of $W$ till first hitting some other element of $W$. 
\end{enumerate}

Matrix $\overline{P}$ is non-negative and, according to 
lemmas \ref{lemma:L02}, \ref{lemma:L04}, \ref{lemma:L06} and definition \eqref{Formula_07}, 
for any $a \in W$ satisfies condition

\beq \label{Formula_08}
	\sum_{b \in W}{ \overline{p}_{a,b}} = 1
\eeq

\begin{lemma} \label{lemma:L07} 
If finite Markov chain $\mathbf{M}$ is irreducible, $W \subseteq Z$ and $|W|>1$, then factor-chain $\mathbf{\overline{M}} = \mathbf{M}/W$  is irreducible too.
\end{lemma}

\begin{proof}
Consider any $a,b \in W$, $a \neq b$. Because $\mathbf{M}$ is irreducible, there exists some path  $x=(x_0,\ldots,x_m)$ that leads (in $\mathbf{M}$) from state $a$  to state $b$  and has positive probability, i.e. 
$p( x_{i-1},x_i ) > 0$ for all $i \leq m$. 

Select in the sequence of states $x$ sub-sequence 
$y=(y_0,\ldots,y_k)=(x_{h(0)},\ldots,x_{h(k)})$
in the following way. Let $H(0)=0$. If $h(j)$ already selected and 
$h(j)<m$, then select  
\beqnn
	h(j+1)=Min\{i|(i>h(j)) \& (x_i \in W) \& (x_i \neq x_{h(j)}) \}
\eeqnn

Evidently, $y$ is subsequence of sequence $x_0,\ldots,x_m$ that contains only elements of $W$, each element in $y$ differs from the previous one, $y_0=x_{h(0)}=a$ and $y_k=y_{h(k)}=b$.

According to definition of $\mathbf{\overline{M}}$, 
corollary  \ref{Cor02} and lemma \ref{lemma:L06}, for any $j\leq k$
\beqnn
	\overline{p}(y_{j-1},y_j)=P(R_{y_{j-1},y_j}^W)=
\overline{p}_{{}_{W \backslash \{ y_{j-1} \},y_{j-1},y_j}} \geq
P_{y_{j-1}}^t(S_{(y_{j-1},* )[+W \backslash \{y_{j-1}\} : y_j]}^t)
\eeqnn 
where $t=h(j)-h(j-1)$. 

On the other hand,
$$	
S_{(y_{j-1},* )[+W \backslash \{y_{j-1}\} : y_j]}^t = 
S_{(x_{h(j-1)},* )[+W \backslash \{x_{h(j-1)}\} : x_{h(j)}]}^t 
$$  
and 
$$
	(x_{h(j-1)},\ldots x_{h(j)}) \in 
	S_{(x_{h(j-1)},* )[+W \backslash \{x_{h(j-1)}\} : x_{h(j)}]}^t
$$ 
Hence, 
\beqnn
	\overline{p}(y_{j-1},y_j) \geq P_{x_{h(j-1)}}^t((x_{h(j-1)},\ldots,x_{h(j)})) =
	\prod_{j<i\leq h(j)}{p(x_{i-1},x_i)}>0
\eeqnn 
and $y=(y_0, \ldots ,y_k)$  is a path in $\mathbf{\overline{M}}$ that leads from $a$ to $b$ and has positive transition probability from each its state to the next one. Proof completed.
\end{proof}

Note that if $p_{z,z}=0$ for all $z \in Z$ and $W=Z$, then chains $\mathbf{M}$ and $\mathbf{\overline{M}}$ coincide. 
From this evident statement follows that transitivity matrix $\overline{P}$ of factor chain could have zero elements outside main diagonal. We could use as a sample $\mathbf{M}/Z$ for any irreducible Markov chain $\mathbf{M}$, having matrix $P$ with the same property. 

We will use dashed letters for all notations, related to factor chain. For example, $\overline{S}$ is the set of all paths in $\mathbf{\overline{M}}$. For $a=b$, according to \eqref{Formula_07}, $\overline{p}_{a,b}=0$, therefore paths with two coinciding adjacent states have zero probability. In following sections we will consider only paths from $\overline{S}^+$:
\beqnn
	\overline{S}^+=\{ x \in \overline{S} | (Len(x)>0) \& (\forall i\leq Len(i))(x_{i-1} \neq x_i) \}
\eeqnn  


\section{Weighted transitions and weighted hitting time}

Hitting time in Markov chain $\mathbf{M}$ defined for hitting some set of states (see \cite{cite_1},p.12), but we will consider hitting time only for sets that contain one state. For any $a,b \in Z$, \emph{hitting time from $a$  to $b$} is random value $\tau_{a,b}$ - minimal index of step when $\mathbf{M}$ achieved state $b$, under condition that initial state is $a$. 

Subject of our interest is the mean hitting time, i.e. function
\beqnn
	H(a,b)=E(\tau_{a,b})=\sum_{i=0}^{\infty}{i\cdot Probability(\tau_{a,b}=i)}
\eeqnn
where $E$ denotes mathematical expectation. 

Note that some authors (e.g. \cite{cite_2},p.29) use a little different definition of hitting time  
$\tau_{a,b}^{'}$: minimal \emph{positive} index of step when $\mathbf{M}$ achieved state $b$, under condition that initial state is $a$. Evidently, $\tau_{a,b}^{'}=\tau_{a,b}$ for all if $a\neq b$. While $\tau_{a,a}=0$ for all $a$, value of $\tau_{a,a}^{'}$ could be positive.
\\

For any set of paths $X$ and any $a,b \in Z$  denote by $X_{<a,b>}$  the set of paths from $X$ that leads from $a$ to $b$, and contains $b$ only as its end state:
\beqnn
	X_{<a,b>}=\{ x \in X | (x \in S_{(a,b)}) \& (\forall i<Len(x))(x_i \neq b) \}
\eeqnn
Evidently, $X_{<a,b>}$ is the set of all $\{ a,b \}$-arrows from $a$ to $b$ that belongs to $X$, i.e. 
\beqnn
	X_{<a,b>}= X \cap R_{a,b}^{\{a,b \}}
\eeqnn
If initial state of Markov chain is $a$, then the sequence of states till first hitting $b$ will follow some random path from $S_{<a,b>}$, and $\tau_{a,b}$ is the length of this path, 
\beqnn
	Probability(\tau_{a,b}=i)= P_a^i(S_{<a,b>}^i)
\eeqnn
and 
\beqnn
	H(a,b)=
	E(\tau_{a,b})=
	\sum_{x \in S_{<a,b>}}
	{Len(x) \cdot P_a^{Len(x)} (x) } 
\eeqnn

For proving some properties of function $H(a,b)$  we need to expand considered notion and introduce weighted hitting time.

Select some positive matrix $V=(v_{i,j})_{i,j=1,\ldots,n}$ ; $(\forall i,j=1,\ldots,n)(v_{i,j}>0)$.
Similar to matrix $P$, notations $v(z_j,z_j)$ and $v_{i,j}$ considered as equivalent.  For any path $x \in S$ define \emph{weight of path $x$ with weight matrix $V$} as  
\beqnn
	Weight(V,x)=\sum_{i=1}^{Len(x)}{v(x_{i-1},x_i)}
\eeqnn

Define $\tau_{V,a,b}$ as weight of random sequence of states of Markov chain with initial state $a$ till first hitting the state $b$. Mathematical expectation of $\tau_{V,a,b}$ is \emph{mean weighted hitting time from $a$  to $b$ with weight matrix $V$}: 
\beq \label{Formula_09}
	H(V,a,b) = \sum_{x \in S_{<a,b>}}{ Weight(V,x) \cdot P_a^{Len(x)} (x)}
\eeq

Let $E$ is a trivial weight matrix, with all elements equal to 1. Evidently,
$$Weight(E,x)=Len(x)$$ and $$H(E,a,b)=H(a,b)$$

\begin{lemma} \label{lemma:L08} 
If finite Markov chain is irreducible, then for any positive weight matrix $V$ and states $a,b \in Z$ mean weighted hitting time from $a$ to $b$ is finite, i.e. the sum 
\beqnn
	H(V,a,b) = \sum_{x \in S_{<a,b>}}{ Weight(V,x) \cdot P_a^{Len(x)} (x)}
\eeqnn 
is convergent.
\end{lemma}

\begin{proof}
Order addends of sum according to length of paths:
$$
	H(V,a,b) = \sum_{x \in S_{<a,b>}}{ Weight(V,x)  P_a^{Len(x)} (x)}=
\sum_{m=0}^{\infty}{\sum_{x \in S_{<a,b>}^m} { Weight(V,x)  P_a^{Len(x)} (x)  }   }
$$  
Let $v_{max}=max\{ v_{i,j} | i,j=1,\ldots , n \} $. Note that for $x \in S_{<a,b>}^m $ holds 
$$	
Weight(V,x) \leq v_{max} \cdot Len(x) 
$$
Evidently, 
$
	S_{<a,b>}^m \subseteq Ext^m( S_{(a,*)[-b]}^{m-1} )  
$. 
Using lemma \ref{lemma:L01}, obtain 
$$
	P_a^m(S_{<a,b>}^m) \leq P_a^m(Ext^m( S_{(a,*)[-b]}^{m-1} )) = P_a^{m-1}(S_{(a,*)[-b]}^{m-1})
$$
and 
$$
	\sum_{x \in S_{<a,b>}^m }{ Weight(V,x)  P_a^{Len(x)} (x) } \leq
v_{max}  m	 \sum_{x \in S_{<a,b>}^m }{P_a^m(x)}=
$$
$$
	v_{max} m	 P_a^m(S_{<a,b>}^m ) 
	\leq  
	v_{max}  m	 P_a^{m-1}(S_{(a,*)[-b]}^{m-1} )	
$$

According to 
corollary \ref{Cor01}, there are $t_0$, $\alpha >0$, $1>\beta >0$ such that for any $t \geq t_0$  
\beqnn
	P_a^t(S_{(a,*)[-b]}^t) \leq \beta^{\alpha t} 
\eeqnn
Hence, for any $t \geq t_0$  
\beqnn
	\sum_{x \in S_{<a,b>}^m}{ Weight(V,x)  P_a^{Len(x)} (x) } 
\leq v_{max}  m  \beta^{\alpha (m-1)}
\eeqnn
Sum 
\beqnn
	\sum_{m=t_0}^{\infty}{m\beta^{\alpha (m-1)} }
\eeqnn
is convergent, and this proves the lemma.   
\end{proof}

\section{Direct sum and concatenation on sets of paths}

Concatenation operation $\otimes$, considered above as operation on paths, could be considered as partial operation on sets of paths. If $a \in Z$, $X \subseteq S_{(*,a)}^*$  and $Y \subseteq S_{(a,*)^*}$, define \emph{concatenation of sets $X$ and $Y$} as $X \otimes Y = \{ x \otimes y | (x \in X)\&(y \in Y) \} $.

Besides, we need another partial operation on sets of paths: \emph{ direct sum $\oplus$}, or union of disjoint sets.
Formal definitions of these partial operations are as follows:

\beq \label{Formula_10}
	W=X \oplus Y \iff (X\cap Y=\oslash) \& (Z=X \cup Y)
\eeq

\begin{align}  \label{Formula_11}
W=X & \otimes Y \iff   \\
		& (\exists a\in Z)( (X \subseteq S_{(*,a)}^*) \& (Y \subseteq S_{(a,*)}^*) \& 
(Z=\{x\otimes y | (x \in X) \& (y \in Y) \})  )													 \nonumber 
\end{align}

For any positive weight matrix $V$ and any path $x$ define function 
\beq \label{Formula_12}
	H^V(x)=Weight(V,x) \cdot P_a^{Len(x)}=
 \Big(  \sum_{i=1}^{Len(x)}  {v(x_{i-1},x_i)}  \Big)
\cdot  \prod_{i=1}^{Len(x)}{p(x_{i-1},x_i)} 
\eeq
	 						
Expand function $H^V(x)$ on sets of paths: for any set of paths $X \subseteq S$ define 

\beq \label{Formula_13}
	H^V(X)=\sum_{x \in X}{H^V(x)}
\eeq

Using definition \eqref{Formula_02_5} of function $P_a^m(x)$, define function $P$ on any path $x$  
\beq \label{Formula_14}
	  P(x)=P_{Begin(x)}^{Len(x)} (x) = \prod_{i=1}^m{p(x_{i-1},x_i)}
\eeq
and expand function $P$ on sets of paths:
\beq \label{Formula_15}
	P(X)=\sum_{x \in X}{P(x)}
\eeq

We consider partial function $P$ (partial function $H^V$) defined only on sets $X \subseteq S$, for which sum, used in definition, is convergent.

Note that 
\beq \label{Formula_15_5}
	H(V,a,b)=H^V(S_{<a,b>})
\eeq

\begin{lemma} \label{lemma:L09} 
If finite Markov chain is irreducible, then for any positive weight matrix $V$ partial operations 
$\oplus$, $\otimes$ and partial functions $P(\cdot)$,$H^V(\cdot)$ satisfy the following properties:

\begin{tabular}{@{~~~}l@{~~~}l@{}}
(a) & $ X \oplus Y = Y \oplus X $																				\\
(b) & $(X \oplus Y)\oplus W = X \oplus (Y \oplus W)$										\\				
(c) & $(X \otimes Y)\otimes W = X \otimes (Y \otimes W)$								\\
(d) & $ X \otimes (Y \oplus W) = (X \otimes Y) \oplus (X \otimes W) $		\\
(e) & $ (X \oplus Y) \otimes W = (X \otimes W) \oplus (Y \otimes W) $		\\
(f) & $ P(X \oplus Y) = P(X)+P(Y) $																			\\
(g) & $ H^V(X \oplus Y) = H^V(X)+H^V(Y) $																\\			
(h) & $ P(X \otimes Y) = P(X) \cdot P(Y) $															\\		
(i) & $ H^V(X \otimes Y) = H^V(X) \cdot P(Y) + P(X) \cdot H^V(Y) $			\\
\end{tabular}
\end{lemma}

As usual, equality of partial functions means that 
\\* - Left and right parts of equality defined or undefined simultaneously and
\\* - If left and right parts are both defined, their values coincide. 

\begin{proof}	
Almost all equalities immediately follow from definition of considered partial operations and partial functions. Verify statements (h) and (i). Let $a\in Z$, $X \subseteq S_{(*,a)}^*$  and $Y \subseteq S_{(a,*)}^*$.

Check (h): 
\begin{align}
P(X \otimes Y) 
	& = \sum_{w \in X\otimes Y}{P(w)} = \sum_{x\in X}{\sum_{y\in Y}{P(x\otimes y)}}= 
		\sum_{x\in X}{\sum_{y\in Y}{ \Big(
		\prod_{i=1}^{Len(x)}{p(x_{i-1},x_i)} 
		\prod_{j=1}^{Len(y)}{p(y_{j-1},y_j)}\Big) }}=   \nonumber \\
	& \sum_{x\in X}{\Big(\prod_{i=1}^{Len(x)}{p(x_{i-1},x_i)}   
		\sum_{y\in Y}{\prod_{j=1}^{Len(y)}{p(y_{j-1},y_j)}}\Big)}=
		\sum_{x\in X}{\Big(\prod_{i=1}^{Len(x)}{p(x_{i-1},x_i)}   
		\sum_{y\in Y}{ P(y)  }\Big)}= \nonumber \\ 
	& P(Y) \sum_{x\in X}{\Big(\prod_{i=1}^{Len(x)}{p(x_{i-1},x_i)}\Big) } =
		P(Y) \sum_{x\in X}{P(X)} = 
		P(Y)\cdot P(X)										\nonumber 
\end{align}

Check (i):
\begin{align}
& 
H^V(X\otimes Y) = \sum_{w \in X \otimes Y}{H^V(w)} = \sum_{x \in X}{\sum_{y \in Y}{H^V(x\otimes y)}}=  \nonumber \\
&
\sum_{x \in X}{\sum_{y \in Y}{
\Big(
\Big(
\sum_{i=1}^{Len(x)}{v(x_{i-1},x_i)}+
\sum_{j=1}^{Len(y)}{v(y_{j-1},y_j)}
\Big)
\prod_{i=1}^{Len(x)}{p(x_{i-1},x_i)} 
\prod_{j=1}^{Len(y)}{p(y_{j-1},y_j)}
\Big)
}}=    \nonumber 
\end{align}

\begin{align}
&
\sum_{x \in X}{\sum_{y \in Y}{
\Big(
\Big(
\sum_{i=1}^{Len(x)}{v(x_{i-1},x_i)}
\Big)
\prod_{i=1}^{Len(x)}{p(x_{i-1},x_i)} 
\prod_{j=1}^{Len(y)}{p(y_{j-1},y_j)}
\Big)
}}+  \nonumber \\
&
\sum_{x \in X}{\sum_{y \in Y}{
\Big(
\Big(
\sum_{j=1}^{Len(y)}{v(y_{j-1},y_j)}
\Big)
\prod_{i=1}^{Len(x)}{p(x_{i-1},x_i)} 
\prod_{j=1}^{Len(y)}{p(y_{j-1},y_j)}
\Big)
}}=      \nonumber 
\end{align}

\begin{align}
&
\sum_{y \in Y}{\Big(
\prod_{j=1}^{Len(y)}{p(y_{j-1},y_j)}
\Big)} 
\sum_{x \in X}{\Big(
\Big(
\sum_{i=1}^{Len(x)}{v(x_{i-1},x_i)}
\Big)
\prod_{i=1}^{Len(x)}{p(x_{i-1},x_i)}
\Big)}+    				\nonumber \\
&
\sum_{x \in X}{\Big(
\prod_{i=1}^{Len(x)}{p(x_{i-1},x_i)}
\Big)} 
\sum_{y \in Y}{\Big(
\Big(
\sum_{j=1}^{Len(y)}{v(y_{j-1},y_j)}
\Big)
\prod_{j=1}^{Len(y)}{p(y_{j-1},y_j)}
\Big)}=     \nonumber 
\end{align}

\begin{align}
&
\sum_{y \in Y}{\Big(
\prod_{j=1}^{Len(y)}{p(y_{j-1},y_j)}
\Big)} \cdot
H^V(X) + 
\sum_{x \in X}{\Big(
\prod_{i=1}^{Len(x)}{p(x_{i-1},x_i)}
\Big)} \cdot
H^V(Y)=     \nonumber \\
&
P(Y)\cdot H^V(X) + P(X) \cdot H^V(Y) \nonumber 
\end{align}

Proof completed.
\end{proof}

Partial operation of direct sum could be expanded on several disjoint sets:
\beqnn
	Y=\bigoplus_{i=1}^m{X_i} 
	\iff 
	\Big(  Y=\bigcup_{i=1}^m{X_i} \Big) \&
	(\forall i,j \leq m)(i\neq j \Rightarrow  X_i\cup X_j = \oslash)
\eeqnn
According to associate property of concatenation (lemma \ref{lemma:L09},c), concatenation of several sets
$$ Y=\bigotimes_{i=1}^m{X_i} = X_1 \otimes (X_2 \otimes (\ldots \otimes X_m)\ldots ) $$  
does not depends on order of parentheses. Similar to concatenation of two sets, 
\begin{align}
Y=\bigotimes_{i=1}^m{X_i} & \iff     \\
& 
(\forall i=1,\ldots ,m-1)(\exists a_i \in Z) (( X_i \subseteq S_{(*,a_i)}^*) \& ( X_{i+1} \subseteq S_{(a_i,*)}^*)) \& 
\nonumber \\
&
(Z = \{ x_1 \otimes x_2 \ldots \otimes x_m | (x_1 \in X_1) \& \ldots \& (x_m \in X_m) \})
\nonumber
\end{align}

\begin{lemma} \label{lemma:L10} 
If finite Markov chain is irreducible, then for any positive weight matrix $V$ and any $m \geq 2$  the following properties hold:\\
(a) $$ P \Big( \bigoplus_{i=1}^m{X_i} \Big) = \sum_{i=1}^m{P(X_i)} $$  
(b) $$ H^V \Big( \bigoplus_{i=1}^m{X_i} \Big) = \sum_{i=1}^m{H^V(X_i)} $$  				
(c) $$ P \Big( \bigotimes_{i=1}^m{X_i} \Big) = \prod_{i=1}^m{P(X_i)} $$  
(d) $$ H^V \Big( \bigotimes_{i=1}^m{X_i} \Big) = 
		\sum_{i=1}^m{\Big( H^V(X_i)\cdot \prod_{{}_{j\neq i}^{j=1}}^m{P(X_j)}   \Big) } $$	
\end{lemma}

\begin{proof}
Statements (a),(b),(c) immediately follow from lemma \ref{lemma:L09} ((f), (g), (h) respectively).  
Prove (d) by induction by $m$. For $m=2$ statement coincides with lemma \ref{lemma:L09}(i). 
If (d) true for some value of $m$, then it is true for $m+1$ too, because
$$
	H^V \Big( \bigotimes_{i=1}^{m+1}{X_i} \Big) =
H^V \Big(\Big( \bigotimes_{i=1}^m{X_i} \Big) \otimes X_{m+1} \Big) =
H^V \Big( \bigotimes_{i=1}^m{X_i} \Big)  P(X_{m+1})+   
H^V(X_{m+1})  P\Big(  \bigotimes_{j=1}^m{X_j}  \Big)=
$$
$$
	\sum_{i=1}^m{\Big(  H^V(X_i)  \prod_{{}_{j\neq i}^{j=1}}{P(X_j)} \Big) }  P(X_{m+1})+
H^V(X_{m+1}) \prod_{j=1}^m{P(X_j)}= 	
$$ 
$$
	\sum_{i=1}^m {\Big(  H^V(X_i) \prod_{{}_{j\neq i}^{j=1}}{P(X_j)} \Big) }+
H^V(X_{m+1})  \prod_{{}_{j\neq m+1}^{j=1}}^{m+1}{P(X_j)}=	
	\sum_{i=1}^{m+1}{\Big( H^V(X_i) \prod_{{}_{j\neq i}^{j=1}}^{m+1}{P(X_j)} \Big)}
$$ 
Proof completed.
\end{proof}

\section{Relation between chain and its factor chain}

Consider finite irreducible Markov chain $\mathbf{M}$ with set of states $Z$, some subset $W\subseteq Z$, $|W|>1$ and factor chain $\mathbf{\overline{M}}= \mathbf{M}/W $. We will use dashed letters for all notations, related to factor chain. 

Let $S_W$ is the set of all paths $x$ in $\mathbf{M}$ that 
\\* - Have non-zero length,
\\* - Connect states from $W$ and
\\* - Last two elements of $x$ that belong to $W$ are different. 

Formally,
$$
	S_W=\{x | (x\in S)\&(Len(x)>0) \&(Begin(x)\in W)\&
$$
$$	
(End(x)\in W) \&(x_{{}_{  max\{ i<Len(x)|x_i \in W \} } } \neq End(x)) \}
$$
 
Note that $S_W$ is closed under partial operation $\otimes$, i.e.
\beqnn
	(v=x\otimes y)\&(x,y\in S_W) \Longrightarrow (v \in S_W)
\eeqnn
 
Define mapping $ \varphi: S_W \longrightarrow \overline{S}$ of $S_W$ in the set $\overline{S}$ of 
all paths in $\mathbf{\overline{M}}$.

Let $x=(x_0,\ldots,x_m) \in S_W$, $a=Begin(x)\in W$ and $b=End(x)\in W$.
Like in the proof of lemma \ref{lemma:L07}, select in the sequence of states $x$ subsequence 
$y=(y_0,\ldots ,y_k)=(x_{h(0)},\ldots ,x_{h(k)})$ 
as follows: 
\\
\\* Let $h(0)=0$.
\\* If $h(j)$ already selected, $h(j)<m$ and 
$$
	(\exists i)((i>h(j)) \&  (x_i \in W) \& (x_i \neq x_{h(j)})  )
$$ 
then 
$$
	  h(j+1) = Min \{ i | (i>h(j)) \&  (x_i \in W) \& (x_i \neq x_{h(j)})   \} 
$$ 
\\
Sequence $y$  satisfies properties: 
\\* - $y$ contains only elements of $W$;
\\* -	Each element in $y$ differs from the previous one;
\\*	- $y_0=x_{h(0)}=a$;
\\* - $y_k=x_{h(k)}=b$.  
\\
Define $\varphi (x) = y$. Here $y\subseteq W$  considered as path in $\mathbf{\overline{M}}$.

Immediately from definition follow that the range of $\varphi$ is the set $\overline{S}^+$ 
of all paths from $\overline{S}$, having non-zero length and different adjacent states: 
\beqnn
	Range(\varphi)=\overline{S}^+=
	\{\overline{x}\in \overline{S} | (Len( \overline{x})>0) \& (\forall i \leq Len(\overline{x}))
	(\overline{x}_{i-1} \neq \overline{x}_i) \}
\eeqnn
 
It is also evident that  $\varphi$ preserves concatenation: 
\beqnn
	(\forall v,x,y \in S_W)(v=x\otimes y \Longrightarrow \varphi(v) = \varphi (x) \otimes \varphi (y))
\eeqnn
	
Hence, $\varphi$  is homomorphism of algebraic system $<S_W, \otimes>$ on system 
$<\overline{S}^+, \otimes >$.
Using homomorphism $\varphi$, we could define on $S_W$ equivalence relation $=_{\varphi}$		 
\beqnn
	(x' =_{\varphi} x'') \iff (\varphi(x')=\varphi(x'')) \qquad x',x''\in S_W
\eeqnn
And equivalence classes
\beqnn
	 [x]_{\varphi}=\{x'\in S_W | x'=_{\varphi}x \}       \qquad    x\in S_W
\eeqnn

Consider again subsequence $y=(y_0,\ldots ,y_k)=(x_{h(0)},\ldots ,x_{h(k)})$, constructed in definition of $\varphi$.
According to definition of function $h$, for all $i=1,\ldots ,k$ segment $(x_{h(i-1)},\ldots ,x_{h(i)})$ of 
sequence $x$ satisfy properties:
\\* -	Its begin state is $x_{h(i-1)}\in W$ ;
\\* -	Its end state is $x_{h(i)}\in W$ ;
\\* -	$x_{h(i-1)} \neq x_{h(i)}$ ;
\\* -	All elements of segment, besides end one, belong to $(Z \backslash W)\cup \{ x_{h(i-1)} \}$ .

This means that $(x_{h(i-1)},\ldots ,x_{h(i)})$ is $W$-arrow from $x_{h(i-1)}$ to $x_{h(i)}$:
\beqnn
	(x_{h(i-1)},\ldots ,x_{h(i)}) \in R_{x_{h(i-1},x_{h(i)} }^W
\eeqnn
and 
\beqnn
	x \in \bigotimes_{i=1}^k  { R_{x_{h(i-1)},x_{h(i)}}^W }
\eeqnn
Note that 
$	\varphi((x_{h(i-1)},\ldots ,x_{h(i)})) = (x_{h(i-1)},x_{h(i)}) $ 
and 
$R_{x_{h(i-1)},x_{h(i)}}^W$ 
is equivalence class of relation $=_{\varphi}$ :
\beqnn
	 R_{x_{h(i-1)},x_{h(i)}}^W =
[(x_{h(i-1)},\ldots ,x_{h(i)})]_{\varphi} =
\varphi^{-1}((x_{h(i-1)},x_{h(i)}))
\eeqnn 
Here $(x_{h(i-1)},x_{h(i)})$ in the rightmost expression considered as a path of length 1 in $\mathbf{\overline{M}}$.
\\

Now we will fix weight matrices for chains $\mathbf{M}$ and $\mathbf{\overline{M}}$. For $\mathbf{M}$ we will use trivial weight matrix $E$, all elements of $E$ equal to 1. Weight matrix 
$\overline{V}=(\overline{v}_{a,b})_{a,b \in W}$ for $\mathbf{\overline{M}}$ 
define as follows: for all couples of states $a,b \in W$ in $\mathbf{\overline{M}}$
\beq \label{Formula_15_6} 
	\overline{v}_{a,b}=
\begin{cases}
H^E(R_{a,b}^W)/P(R_{a,b}^W) & \text{if } \quad (a\neq b)\&(P(R_{a,b}^W)>0)\\
1 & \text{otherwise} 
\end{cases}
\eeq

This means that for all $a,b\in W$, $a\neq b$, if $P(R_{a,b}^W)>0$ , then   
\beqnn
	\overline{v}_{a,b}=\Big(   \sum_{x\in R_{a,b}^W}{Len(x) P_a^{Len(x)}(x)}  \Big) / 
\Big(  \sum_{x\in R_{a,b}^W}{P_a^{Len(x)}(x)}  \Big)	
\eeqnn 
 
Note that if $a,b\in W$, $a\neq b$  and $P(R_{a,b}^W)=0$, then $P_a^{Len(x)}(x)=0$  for all $x\in R_{a,b}^W$ 
and hence $H^E(R_{a,b}^W)=0$  too.

Let $H^E(\cdot)$  is the function, defined on paths and set of paths in $\mathbf{M}$ according to \eqref{Formula_12} and \eqref{Formula_13} with usage of trivial weight matrix $E$, and $\overline{H}^{\overline{V}}(\cdot)$  is the similar function in $\mathbf{\overline{M}}$, using weight matrix $\overline{V}$ .

\begin{lemma} \label{lemma:L11} 
If finite Markov chain $\mathbf{M}$ is irreducible, $\mathbf{\overline{M}}=\mathbf{M}/W$
is factor chain, \mbox{$|W|>1$}, and functions $\varphi$, $H^E(\cdot)$,
$\overline{H}^{\overline{V}}(\cdot)$ are as defined above, then\\
(a)	For all $a,b\in W$,$a\neq b$ 
$$\overline{H}^{\overline{V}}((a,b))=H^E(\varphi^{-1}((a,b)))$$ 
(b) For all $a,b,c\in W$, $a\neq b$, $b\neq c$
$$\overline{H}^{\overline{V}}((a,b)\otimes (b,c))= H^E(  \varphi^{-1}((a,b)) \otimes \varphi^{-1}((b,c)) )$$
(c)	For any path $y=(y_0,\ldots , y_k)\in \overline{S}$, $(\forall i \leq k)(y_{i-1} \neq y_i)$
$$  
\overline{H}^{\overline{V}}(y)=H^E \Big(   
\bigotimes_{i=1}^k{ \varphi^{-1}((y_{i-1},y_i))}
\Big)
$$ 
\end{lemma}

\begin{proof}
(a)  For any $a,b \in W$, $a\neq b$ if $P(R_{a,b}^W)>0$, then 
$$
\overline{H}^{\overline{V}}((a,b))=
\overline{v}_{a,b}\cdot \overline{p}_{a,b}=
(H^E(R_{a,b}^W)/P(R_{a,b}^W))\cdot P(R_{a,b}^W) =
H^E(R_{a,b}^W)=
H^E(\varphi^{-1}((a,b)))
$$
If $P(R_{a,b}^W)=0$, then 
$\overline{H}^{\overline{V}}((a,b))= \overline{v}_{a,b}\cdot \overline{p}_{a,b}= 0$ 
and $H^E(\varphi^{-1}((a,b))) =H^E(R_{a,b}^W)=0$.
\\*

(b) Note that for any $a,b \in W$, $a\neq b$ 
$$
\overline{p}_{a,b}=P(R_{a,b}^W)=P(\varphi^{-1}((a,b)))
$$

Let $a,b,c\in W$, $a\neq b$, $b\neq c$. Applying lemma \ref{lemma:L09}(i) to concatenation of one-element sets 
of paths $(a,b)\otimes (b,c)$ in $\mathbf{\overline{M}}$, receive  
\begin{align} 
\overline{H}^{\overline{V}}((a,b)\otimes (b,c)) & =
\overline{H}^{\overline{V}}((a,b)) \cdot \overline{P}((a,b))+
\overline{P}((a,b)) \cdot \overline{H}^{\overline{V}}((b,c))=  \nonumber \\
&
H^E(\varphi^{-1}((a,b))) \cdot \overline{p}_{b,c}+
\overline{p}_{a,b} \cdot H^E(\varphi^{-1}((b,c)))=  \\
&
H^E(\varphi^{-1}((a,b))) \cdot P(\varphi^{-1}((b,c))) +
P(\varphi^{-1}((a,b))) \cdot H^E(\varphi^{-1}((b,c)))  \nonumber
\end{align}
Now (b) follows from lemma \ref{lemma:L09}(i) for sets of paths $\varphi^{-1}((a,b))$  
and $\varphi^{-1}((b,c))$ in $\mathbf{M}$.  
\\*

(c) Similarly to proof of (b), for any $y=(y_0,\ldots ,y_k)\in \overline{S}$,
$(\forall i\leq k)(y_{i-1}\neq y_i)$ use lemma \ref{lemma:L10}(d): 
$$
\overline{H}^{\overline{V}}(y)=
\overline{H}^{\overline{V}}\Big( \bigotimes_{i=1}^k{(( y_{i-1},y_i ))} \Big)=
\sum_{i=1}^k{\Big(
	  \overline{H}^{\overline{V}}((y_{i-1},y_i)) 
	  \prod_{{}_{j\neq i}^{j=1}}^m{ \overline{P}((y_{j-1},y_j))  } 
\Big)}=
$$
 
$$
\sum_{i=1}^k{\Big(
	  H^E( \varphi^{-1}(y_{i-1},y_i)) 
	  \prod_{{}_{j\neq i}^{j=1}}^m{ P(\varphi^{-1}(y_{j-1},y_j))} 
\Big)}=
H^E\Big( 
	\bigotimes_{i=1}^k{\varphi^{-1}(y_{i-1},y_i)}
\Big)
$$ 
Proof completed.
\end{proof}

Recall that if weight matrix $V$ is given, weighted hitting time from $a$ to $b$ is 
\beqnn
	H(V,a,b)=\sum_{x\in S_{<a,b>}}{H^V(x)}
\eeqnn
where $S_{<a,b>}=R_{a,b}^{\{ a,b\}}$ is the set of paths that leads from $a$ to $b$ and contain $b$ only as end state.

\begin{lemma} \label{lemma:L12} 
If finite Markov chain $\mathbf{M}$ is irreducible, $\mathbf{\overline{M}}=\mathbf{M}/W$
is factor chain, \mbox{$|W|>1$}, $a,b \in W$ and $a\neq b$, then 
\beqnn
	\overline{H}(\overline{V},a,b)=H(E,a,b)
\eeqnn
\end{lemma}
Here  $\overline{H}(\overline{V},a,b)$ is weighted mean hitting time in factor chain $\mathbf{\overline{M}}$, calculated using matrix $\overline{V}$, defined above \eqref{Formula_15_6}.

\begin{proof}
In equality
\beqnn
	\overline{H}(\overline{V},a,b)=\sum_{x\in \overline{S}_{<a,b>}}{\overline{H}^{\overline{V}} (x) }
\eeqnn
drop from the sum addends for those paths $x$ that have coinciding adjacent states. This does not affect the sum, because these paths have zero probability.  Hence, instead of $\overline{S}_{<a,b>}$ we could use 
\begin{align}
& \overline{S}_{<a,b>}^+ = \overline{S}_{<a,b>} \cap \overline{S}^+ = \\
&
\{x\in \overline{S} |   
(Len(x)>0) \& (Begin(x)=a) \& (End(x)=b) \& 
(\forall i\leq Len(i))(x_{i-1} \notin \{x_i,b \})
\}    \nonumber
\end{align}
According to lemma \ref{lemma:L11}(c),
\beq \label{Formula_16}
	\overline{H}(\overline{V},a,b) = 
\sum_{x\in \overline{S}_{<a,b>}^+}{\overline{H}^{\overline{V}}(x)}= 		
\sum_{x\in \overline{S}_{<a,b>}^+}{H^E \Big( 
	\bigotimes_{i=1}^{Len(x)}{ \varphi^{-1}((x_{i-1},x_i )) }
\Big)}
\eeq

Note that if $x',x'' \in \overline{S}_{<a,b>}^+$ and $x' \neq x''$  , then sets 
$$ 
\bigotimes_{i=1}^{Len(x')}{ \varphi^{-1} ((x_{i-1}^{'},x_i^{'} )) }
\qquad \text{ and } \qquad
\bigotimes_{i=1}^{Len(x'')}{ \varphi^{-1}((x_{i-1}^{''},x_i^{''} )) } 
$$
are disjoint. Actually, because  $\varphi$ preserves concatenation, any path $u$ that belongs to intersection of these two sets must satisfy inconsistent system of equalities: $\varphi(u)=x'$  and $\varphi(u)=x''$.

Applying to  \eqref{Formula_16} lemma \ref{lemma:L10}(b), obtain
\beqnn
	\overline{H}(\overline{V},a,b) = H^E\Big( 
	\bigoplus_{x\in \overline{S}_{<a,b>}^+}{ \bigotimes_{i=1}^{Len(x)}{ \varphi^{-1}((x_{i-1},x_i)) }}
	\Big)
\eeqnn
For completion proof it is sufficient verify equality 
\beqnn
	S_{<a,b>}=\bigoplus_{x\in \overline{S}_{<a,b>}^+}{ \bigotimes_{i=1}^{Len(x)}{ \varphi^{-1}((x_{i-1},x_i)) }}
\eeqnn
or
\beqnn
	S_{<a,b>}=\bigoplus_{x\in \overline{S}_{<a,b>}^+}{ \bigotimes_{i=1}^{Len(x)}{ R_{x_{i-1},x_i}^W }}
\eeqnn
If $x=(x_0,\ldots ,x_k)\in \overline{S}_{<a,b>}^+ $ , then all elements in $x=(x_0,\ldots ,x_k)$  belong to $W$,
$a=x_0$, each element differ from the previous one, and $x_k$ is the only element that coincide with $b$. 
This means that any path from 
$$
\bigotimes_{i=1}^k{R_{x_{i-1},x_i}^W}
$$
leads (in chain $\mathbf{M}$) from $a$ to $b$, and contains $b$ only at its end. Hence, $x \in S_{<a,b>}$.  

On the other hand,if $x=(x_0,\ldots ,x_m)\in S_{<a,b>}$  then, using subsequence
$y=(y_0,\ldots ,y_k)=(x_{h(0)},\ldots , x_{h(k)})$
constructed for definition $\varphi$, receive  
$$
	x=\bigotimes_{i=1}^k{(x_{h(i-1)},\ldots , x_{h(i)})} 
$$  
and 
$$
	(x_{h(i-1)},\ldots , x_{h(i)}) \in R_{x_{h(i-1)},x_{h(i)}}^W 
$$  
Subsequence $y$ was selected in such a way that $(y_0,\ldots ,y_k) \in \overline{S}_{<a,b>}^+$, therefore

\beqnn
	x\in \bigotimes_{i=1}^k{R_{y_{i-1},y_i}^W} \subseteq 
  \bigoplus_{x\in \overline{S}_{<a,b>}^+}{ \bigotimes_{i=1}^{Len(x)}{ R_{x_{i-1},x_i}^W }}
\eeqnn
Proof completed.
\end{proof}

\section{The triangle inequality for mean hitting time}

\begin{theorem} \label{theorem:T01}
In finite irreducible Markov chain, mean hitting time $H$ satisfies triangle inequality 
\beq \label{Formula_17}
	H(a,c) \leq H(a,b)+H(b,c)
\eeq
for any states $a,b,c$.
\end{theorem}

\begin{proof} 

If at least two of states $a,b,c$ coincide, the statement is trivial. Consider non-trivial case of mutually different $a,b,c$.

Let $W=\{a,b,c\}$ . Consider factor chain $\mathbf{\overline{M}}=\mathbf{M}/W$ with weight matrix $\overline{V}$, defined above \eqref{Formula_15_6}.  According to lemma \ref{lemma:L12}, 
$H(a,b)=H(E,a,b)=\overline{H}(\overline{V},a,b)$,  
$H(b,c)=\overline{H}(\overline{V},b,c)$
and $H(a,c)=\overline{H}(\overline{V},a,c)$.
Hence, instead of \eqref{Formula_17} we could prove equivalent inequality 
\beq \label{Formula_18}
	\overline{H}(\overline{V},a,c) \leq 
	\overline{H}(\overline{V},a,b)+
	\overline{H}(\overline{V},b,c)
\eeq
in very simple factor chain $\mathbf{\overline{M}}$: it has only three states $a,b,c$, and 
$\overline{p}_{a,a}=\overline{p}_{b,b}=\overline{p}_{c,c}=0$.
\\

Calculate 	$\overline{H}(\overline{V},a,c)$, $\overline{H}(\overline{V},a,b)$ and 
$\overline{H}(\overline{V},b,c)$.
\beqnn
	\overline{H}(\overline{V},a,c) = 
	\sum_{x \in \overline{S}_{<a,c>}}{ Weight(\overline{V},x) \cdot \overline{P}_a^{Len(x)}(x)}
\eeqnn
Here $\overline{S}_{<a,c>}$ is the set of all paths in
$W=\{ a,b,c \}$   
that leads from $a$ to $c$ and contains $c$ only as its end state.

Note that 
\beq \label{Formula_19}
	\overline{P}(\overline{S}_{<a,c>}) =1
\eeq
Formula \eqref{Formula_19} asserts that chain $\mathbf{\overline{M}}$, being in state $a$, will hit (after one or several steps) state $c$ with probability 1. Actually, 
\beqnn
	\overline{P}(\overline{S}_{<a,c>})=
\sum_{x\in \overline{S}_{<a,c>}}{\overline{P}_a^{Len(x)}(x)}=
\lim_{m \to \infty}{ 
	\sum_{{}_{Len(x)\leq m}^{x\in \overline{S}_{<a,c>}}}
	{\overline{P}_a^{Len(x)}(x)} 
}
\eeqnn
And, according to lemma \ref{lemma:L01},
\begin{align}
\overline{P}(\overline{S}_{<a,c>}) = &
\lim_{m \to \infty}
{
 \sum_{{}_{Len(x)\leq m}^{x\in \overline{S}_{<a,c>}}}
 {\overline{P}_a^m(Ext^m(x))}
}=
\lim_{m \to \infty}{\Big(
  1  - \sum_{  x\in \overline{S}_{(a,*)[-c]}^m  }{\overline{P}_a^m(x)}
\Big)}=  \\
&
1-\lim_{m \to \infty}{\Big(
  \sum_{  x\in \overline{S}_{(a,*)[-c]}^m  }{\overline{P}_a^m(x)}
\Big)}  \nonumber
\end{align}

According to lemma \ref{lemma:L02}
\beqnn
\lim_{m \to \infty}{\Big(
  \sum_{  x\in \overline{S}_{(a,*)[-c]}^m  }{\overline{P}_a^m(x)}
\Big)}=0	
\eeqnn 

and \eqref{Formula_19} proved.
\\

Evidently, in $\mathbf{\overline{M}}$ any path that leads from $a$ till hitting $c$, either coincides with 
$(a,c)$, or starts from $(a,b)$. This means that 

\beq \label{Formula_20}
	\overline{S}_{<a,b>}=(a,c)\oplus \left(  
	(a,b)\otimes \overline{S}_{<b,c>}
\right)
\eeq
According to \eqref{Formula_15_5},
\beq \label{Formula_20_1}
	\overline{H}(\overline{V},a,b)=
	\overline{H}^{\overline{V}}(\overline{S}_{<a,b>})
\eeq
therefore from \eqref{Formula_20} obtain
\beq \label{Formula_20_2}
	\overline{H}(\overline{V},a,c)=
	\overline{H}^{\overline{V}}(\overline{S}_{<a,c>})=
	\overline{H}^{\overline{V}}\left(
		 (a,c)\oplus \left((a,b)\otimes \overline{S}_{<b,c>}\right)
	\right)
\eeq
Applying to \eqref{Formula_20_2} lemma \ref{lemma:L09}(g,i), obtain 
\begin{align}
	\overline{H}(\overline{V},a,c) & =
	\overline{H}^{\overline{V}}((a,c))+
	\overline{H}^{\overline{V}}((a,b)\otimes \overline{S}_{<b,c>})=	\\
&
	\overline{H}^{\overline{V}}((a,c))+
	\overline{H}^{\overline{V}}((a,b)) \cdot \overline{P}(\overline{S}_{<b,c>})+
	\overline{P}((a,b)) \cdot \overline{H}^{\overline{V}}(\overline{S}_{<b,c>})  \nonumber 
\end{align}
Using \eqref{Formula_19} and \eqref{Formula_20_1}, receive 
$$
	\overline{H}(\overline{V},a,c)=
	\overline{H}^{\overline{V}}((a,c))+
	\overline{H}^{\overline{V}}((a,b))+
	\overline{P}((a,b)) \cdot \overline{H}(\overline{V},b,c)	
$$ 
Values of function $\overline{H}^{\overline{V}}$ on paths of length 1 directly expressed according \eqref{Formula_12}:
$$
	\overline{H}^{\overline{V}}((a,c)) = \overline{v}_{a,c}\cdot \overline{p}_{a,c}
$$
$$
	\overline{H}^{\overline{V}}((a,b)) = \overline{v}_{a,b}\cdot \overline{p}_{a,b}
$$
therefore
\beq \label{Formula_21}
	\overline{H}(\overline{V},a,c) =
\overline{v}_{a,c}\cdot \overline{p}_{a,c}+
\overline{p}_{a,b}\cdot (
\overline{v}_{a,b}+\overline{H}(\overline{V},b,c)
)
\eeq
Swapping roles of $a$ and $b$, receive 
\beq \label{Formula_22}
	\overline{H}(\overline{V},b,c)= 
\overline{v}_{b,c}\cdot \overline{p}_{b,c}+
\overline{p}_{b,a} \cdot(
\overline{v}_{b,a}+\overline{H}(\overline{V},a,c)
)	
\eeq
Substitute \eqref{Formula_22} into \eqref{Formula_21} and get 
\beq
	\overline{H}(\overline{V},a,c)= 
\overline{v}_{a,c}\cdot \overline{p}_{a,c}+
\overline{p}_{a,b} \cdot  \left(
\overline{v}_{a,b}+
\overline{v}_{b,c}\cdot \overline{p}_{b,c}+
\overline{p}_{b,a} \cdot
(
\overline{v}_{b,a}+
\overline{H}(\overline{V},a,c)
)
\right)	
\eeq
and
\beq \label{Formula_23}
(1- \overline{p}_{a,b} \cdot \overline{p}_{b,a}) \cdot
\overline{H}(\overline{V},a,c) =
\overline{v}_{a,c}\cdot \overline{p}_{a,c}+
\overline{p}_{a,b} \cdot 
(
\overline{v}_{a,b}+
\overline{v}_{b,c} \cdot \overline{p}_{b,c}+
\overline{p}_{b,a} \cdot \overline{v}_{b,a}
)
\eeq

Recall that 
\beq \label{Formula_24}
\begin{cases}
\overline{p}_{a,b}+\overline{p}_{a,c}=1 \\
\overline{p}_{b,a}+\overline{p}_{b,c}=1 \\
\overline{p}_{c,a}+\overline{p}_{c,b}=1 
\end{cases}
\eeq

Note that  
\beqnn 
1- \overline{p}_{a,b}\cdot \overline{p}_{b,a}>0
\eeqnn
 
Actually, otherwise $\overline{p}_{a,b}= \overline{p}_{b,a}=1$, and according to \eqref{Formula_24},
$\overline{p}_{a,c}=0$  and $\overline{p}_{b,c}=0$. These two equalities contradict to irreducibility of 
$\mathbf{\overline{M}}$ , asserted in lemma \ref{lemma:L07}. 

Hence, from \eqref{Formula_23} follow
\beq \label{Formula_25}
\overline{H}(\overline{V},a,c)=\frac
{
	\overline{v}_{a,c}\cdot \overline{p}_{a,c}+
	\overline{p}_{a,b}\cdot \overline{v}_{a,b}+
	\overline{p}_{a,b}\cdot
	(\overline{v}_{b,c}\cdot \overline{p}_{b,c}+\overline{p}_{b,a}\cdot \overline{v}_{b,a})
}
{
	1-\overline{p}_{a,b}\cdot \overline{p}_{b,a}
}
\eeq
Swapping roles of $b,c$ and $a,b$, receive two similar formulas:
\beq \label{Formula_26}
\overline{H}(\overline{V},a,b)=\frac
{
	\overline{v}_{a,b}\cdot \overline{p}_{a,b}+
	\overline{p}_{a,c}\cdot \overline{v}_{a,c}+
	\overline{p}_{a,c}\cdot
	(\overline{v}_{c,b}\cdot \overline{p}_{c,b}+\overline{p}_{c,a}\cdot \overline{v}_{c,a})
}
{
	1-\overline{p}_{a,c}\cdot \overline{p}_{c,a}
}
\eeq

\beq \label{Formula_27}
\overline{H}(\overline{V},b,c)=\frac
{
	\overline{v}_{b,c}\cdot \overline{p}_{b,c}+
	\overline{p}_{b,a}\cdot \overline{v}_{b,a}+
	\overline{p}_{b,a}\cdot
	(\overline{v}_{a,c}\cdot \overline{p}_{a,c}+\overline{p}_{a,b}\cdot \overline{v}_{a,b})
}
{
	1-\overline{p}_{b,a}\cdot \overline{p}_{a,b}
}
\eeq
 
Substituting these three expressions in \eqref{Formula_18} (and dropping overlines for simpler notations) we receive that \eqref{Formula_18} is equivalent to inequality 

\begin{align} \label{Formula_28} 
& 
	\frac{
	v_{a,c}p_{a,c}+p_{a,b}v_{a,b}+p_{a,b}(  v_{b,c}p_{b,c}+ p_{b,a}v_{b,a})   
	}{1-p_{a,b}p_{b,a}} \leq \nonumber \\
& 
	\frac{
	v_{a,b}p_{a,b}+p_{a,c}v_{a,c}+p_{a,c}(  v_{c,b}p_{c,b}+ p_{c,a}v_{c,a})   
	}{1-p_{a,c}p_{c,a}} \\
& 
	\frac{
	v_{b,c}p_{b,c}+p_{b,a}v_{b,a}+p_{b,a}(  v_{a,c}p_{a,c}+ p_{a,b}v_{a,b})   
	}{1-p_{b,a}p_{a,b}} \nonumber 
\end{align}

Because all denominators are positive, for proving \eqref{Formula_28} it is sufficient to verify that expression $g$ is non-negative:
\begin{align} \label{Formula_29} 
g = & (1-p_{a,c}p_{c,a}) 
		\big(v_{b,c}p_{b,c}+p_{b,a}v_{b,a}+p_{b,a}(v_{a,c}p_{a,c}+ p_{a,b}v_{a,b})  \nonumber \\
		& -v_{a,c}p_{a,c}-p_{a,b}v_{a,b}-p_{a,b}(v_{b,c}p_{b,c}+ p_{b,a}v_{b,a})\big)+ \\
		& (1-p_{a,b}p_{b,a})(v_{a,b}p_{a,b}+p_{a,c}v_{a,c}+p_{a,c}(  v_{c,b}p_{c,b}+ p_{c,a}v_{c,a})) \nonumber
\end{align}

Make a series of transformations:
\begin{align}  
g = & (1-p_{a,c}p_{c,a}) 
		(v_{b,c}p_{b,c}+p_{b,a}v_{b,a}+p_{b,a}v_{a,c}p_{a,c}+ p_{b,a}p_{a,b}v_{a,b}  \nonumber \\
		& -v_{a,c}p_{a,c}-p_{a,b}v_{a,b}-p_{a,b}v_{b,c}p_{b,c}-p_{a,b}p_{b,a}v_{b,a})+ \nonumber\\
		& (1-p_{a,b}p_{b,a})(v_{a,b}p_{a,b}+p_{a,c}v_{a,c}+p_{a,c}v_{c,b}p_{c,b}+p_{a,c}p_{c,a}v_{c,a}) \nonumber
\end{align}

\begin{align}  
g = & (1-p_{a,c}p_{c,a}) \cdot  \nonumber \\
		& \left(
			v_{a,b}(p_{b,a}p_{a,b}-p_{a,b}) + v_{a,c}(p_{b,a}p_{a,c}-p_{a,c}) +
			v_{b,a}(p_{b,a}-p_{a,b}p_{b,a}) + v_{b,c}(p_{b,c}-p_{a,b}p_{b,c})
			\right) +  \nonumber \\
		& (1-p_{a,b}p_{b,a}) \cdot 
		(v_{a,b}p_{a,b}+v_{a,c}p_{a,c} +v_{c,a}p_{a,c}p_{c,a} +v_{c,b}p_{a,c}p_{c,b} ) \nonumber
\end{align}
 
\begin{align}  
g = & (1-p_{a,c}p_{c,a}) \cdot  \nonumber \\
		& \left(
			v_{a,b}p_{a,b}(p_{b,a}-1) + v_{a,c}p_{a,c}(p_{b,a}-1) +
			v_{b,a}p_{b,a}(1-p_{a,b}) + v_{b,c}p_{b,c}(1-p_{a,b})
			\right)  + \nonumber \\
		& (1-p_{a,b}p_{b,a}) \cdot 
		(v_{a,b}p_{a,b}+v_{a,c}p_{a,c} +v_{c,a}p_{a,c}p_{c,a} +v_{c,b}p_{a,c}p_{c,b} ) \nonumber
\end{align}

\begin{align}  
g = & (1-p_{a,c}p_{c,a}) \cdot \left(
			-v_{a,b}p_{a,b}p_{b,c} - v_{a,c}p_{a,c}p_{b,c} + v_{b,a}p_{b,a}p_{a,c} + v_{b,c}p_{b,c}p_{a,c}
			\right)+   \nonumber \\
		& (1-p_{a,b}p_{b,a}) \cdot 
		(v_{a,b}p_{a,b}+v_{a,c}p_{a,c} +v_{c,a}p_{a,c}p_{c,a} +v_{c,b}p_{a,c}p_{c,b} ) \nonumber
\end{align}

Consider $g$ as linear form of weights:
\begin{align}  
g = & v_{a,b}(-(1-p_{a,c}p_{c,a})p_{a,b}p_{b,c} + (1-p_{a,b}p_{b,a})p_{a,b} )+  \nonumber \\
		& v_{a,c}(-(1-p_{a,c}p_{c,a})p_{a,c}p_{b,c} + (1-p_{a,b}p_{b,a})p_{a,c} )+  \nonumber \\
		& v_{b,a}(1-p_{a,c}p_{c,a})p_{b,a}p_{a,c}+																	
		  v_{b,c}(1-p_{a,c}p_{c,a})p_{b,c}p_{a,c}+																	\nonumber \\
		& v_{c,a}(1-p_{a,b}p_{b,a})p_{a,c}p_{c,a}+																	
		  v_{c,b}(1-p_{a,b}p_{b,a})p_{a,c}p_{c,b}																	  \nonumber 
\end{align}
and do some transformations, using (24) \eqref{Formula_24}:
\begin{align}  
g = & v_{a,b}(1-p_{b,c}+p_{a,c}p_{c,a}p_{b,c}-p_{a,b}p_{b,a})p_{a,b}+  \nonumber \\
		& v_{a,c}(1-p_{b,c}+p_{a,c}p_{c,a}p_{b,c}-p_{a,b}p_{b,a})p_{a,c}+  \nonumber \\
		& v_{b,a}(1-p_{a,c}p_{c,a})p_{b,a}p_{a,c}+												 
		  v_{b,c}(1-p_{a,c}p_{c,a})p_{b,c}p_{a,c}+												 \nonumber \\
		& v_{c,a}(1-p_{a,b}p_{b,a})p_{a,c}p_{c,a}+												 
		  v_{c,b}(1-p_{a,b}p_{b,a})p_{a,c}p_{c,b}													 \nonumber 
\end{align}
\begin{align}  
g = & v_{a,b}(p_{b,a}+p_{a,c}p_{c,a}p_{b,c}-p_{a,b}p_{b,a})p_{a,b}+  \nonumber \\
		& v_{a,c}(p_{b,a}+p_{a,c}p_{c,a}p_{b,c}-p_{a,b}p_{b,a})p_{a,c}+  \nonumber \\
		& v_{b,a}(1-p_{a,c}p_{c,a})p_{b,a}p_{a,c}+												 
		  v_{b,c}(1-p_{a,c}p_{c,a})p_{b,c}p_{a,c}+												 \nonumber \\
		& v_{c,a}(1-p_{a,b}p_{b,a})p_{a,c}p_{c,a}+												 
		  v_{c,b}(1-p_{a,b}p_{b,a})p_{a,c}p_{c,b}													 \nonumber 
\end{align}
\begin{align}  
g = & v_{a,b}(p_{b,a}(1-p_{a,b})+p_{a,c}p_{c,a}p_{b,c})p_{a,b}+  \nonumber \\
		& v_{a,c}(p_{b,a}(1-p_{a,b})+p_{a,c}p_{c,a}p_{b,c})p_{a,c}+  \nonumber \\
		& v_{b,a}(1-p_{a,c}p_{c,a})p_{b,a}p_{a,c}+												 
		  v_{b,c}(1-p_{a,c}p_{c,a})p_{b,c}p_{a,c}+												 \nonumber \\
		& v_{c,a}(1-p_{a,b}p_{b,a})p_{a,c}p_{c,a}+												 
		  v_{c,b}(1-p_{a,b}p_{b,a})p_{a,c}p_{c,b}													 \nonumber 
\end{align}
\begin{align}  
g = & v_{a,b}(p_{b,a}p_{a,c}+p_{a,c}p_{c,a}p_{b,c})p_{a,b}+  \nonumber \\
		& v_{a,c}(p_{b,a}p_{a,c}+p_{a,c}p_{c,a}p_{b,c})p_{a,c}+  \nonumber \\
		& v_{b,a}(1-p_{a,c}p_{c,a})p_{b,a}p_{a,c}+												 
		  v_{b,c}(1-p_{a,c}p_{c,a})p_{b,c}p_{a,c}+												 \nonumber \\
		& v_{c,a}(1-p_{a,b}p_{b,a})p_{a,c}p_{c,a}+												 
		  v_{c,b}(1-p_{a,b}p_{b,a})p_{a,c}p_{c,b}													 \nonumber 
\end{align}

Now $g$ is linear form on non-negative weights with non-negative coefficients for each weight. 
Hence, $g\geq 0$ and proof completed. 

\end{proof}

\section{Metric on space of states}

\begin{theorem} \label{theorem:T02}
If $\mathbf{M}$ is finite irreducible Markov chain with state space $Z$, then function
$$
	\rho(a,b)=H(a,b)+H(b,a)
$$		 
is a metric on $Z$.
\end{theorem}
\begin{proof} 
The theorem asserts that for all $a,b,c \in Z$  the following properties hold:

\begin{tabular}{@{~~~}l@{~~~}l@{}}
(a) & $ \rho (a,b) \geq 0 $	\\
(b) & $ \rho (a,b)=0 \Leftrightarrow a=b $	\\				
(c) & $ \rho (a,b)= \rho (b,a) $	\\
(d) & $ \rho (a,c) \leq \rho (a,b)+\rho( b,c) $	\\
\end{tabular}

Properties (a)-(c) are evident from definition of $\rho$. 
Property (d) follows from the theorem \ref{theorem:T01}:
\begin{align}  
\rho(a,c) & = H(a,c)+H(c,a) \leq (H(a,b)+H(b,c))+(H(c,b)+H(b,a)) =    \nonumber \\
		    	& (H(a,b)+H(b,a))+(H(b,c)+H(c,b))=\rho(a,b)+ \rho(b,c)      \nonumber 
\end{align}
\end{proof}
Note that $\rho(a,b)$ has alternative equivalent definition
as mean length of minimal loops that include both $a$ and $b$. 
Set of all such loops is $ (S_{<a,b>}\otimes S_{<b,a>}) \cup (S_{<b,a>}\otimes S_{<a,b>})$ .
\\*

Let $Z=\{ z_1,\ldots ,z_n \}$ is the state space of irreducible Markov chain $\mathbf{M}$ with transition matrix $P=(p_{i,j})_{i,j=1,\ldots,n}$.  
Matrix of mean hitting times $h=(h_{i,j})_{i,j=1,\ldots,n}$ (and hence matrix of distances
$\rho=(\rho_{i,j})_{i,j=1,\ldots,n}$) could be calculated, using theorem for calculation mean hitting time of arbitrary subset $A\subseteq Z$  (\cite{cite_1}, Theorem 1.3.5].

This theorem asserts that vector $(k_i^A)_{i=1,\ldots ,n}$ of mean hitting time (form state $z_i$
to some element of $A$) is the minimal non-negative solution of the system of linear equations 

$$
\begin{cases}
k_i^A=0 & \text{if } z_i\in A  \\
k_i^A=1+\sum_{z_j \notin A}{p_{i,j}k_j^A}  & \text{if } z_i\notin A
\end{cases}
$$

Here ``minimal non-negative solution'' means that
$(\forall i=1,\ldots ,n )(k_i^A  \leq  t_i)$
for any other non-negative solution $(t_i)_{i=1,\ldots ,n}$.

For any fixed value of index $j$, apply this theorem to one-element set $A=\{ z_j\}$. 
and receive that vector $(h_{i,j})_{i=1,\ldots ,n}$ is the minimal
non-negative solution of the system of linear equations
$$
\begin{cases}
h_{j,j}=0 &   \\
h_{i,j}=1+\sum_{k \neq j}{p_{i,k}h_k^j}  & \text{for } i \neq j
\end{cases}
$$
Simplified procedures for calculation mean hitting times proposed by Hunter \cite{cite_3}.


\end{document}